\documentstyle[12pt,spacing]{amsart}
\begin{document}
\begin{center}
{\Large\bf A Hierarchy of Maps Between Compacta} 
\\[20pt]
Paul Bankston\\
Department of Mathematics, Statistics and Computer Science\\ 
Marquette University\\
Milwaukee, WI 53201-1881\\[20pt]

\end{center}
\begin{abstract}
Let {\bf CH} be the class of compacta (i.e., compact Hausdorff spaces),
with {\bf BS} the subclass of Boolean spaces.  For each ordinal
$\alpha$ and pair $\langle \mbox{\bf K, L}\rangle$ of subclasses of {\bf CH},
we define $\mbox{Lev}_{\geq\alpha}(\mbox{\bf K,L})$, the class of maps of
{\it level\/} $\geq \alpha$ from spaces in {\bf K} to spaces in {\bf L},
in such a way that  
$\mbox{Lev}_{\geq\alpha}(\mbox{\bf BS,BS})$ consists of the Stone duals of
Boolean lattice embeddings that preserve all prenex first-order formulas
of quantifier rank $\alpha$.  Maps of level $\geq 0$ are just the continuous
surjections, and the maps of level $\geq 1$ are the co-existential maps 
introduced
in \cite{Ban9}.  Co-elementary maps are of level $\geq \alpha$ for all
ordinals $\alpha$; of course in the Boolean context, the co-elementary
maps coincide with the maps of level $\geq \omega$. The results of this paper
include: 
$(i)$ every map of level $\geq \omega$ is co-elementary; 
$(ii)$ the limit maps of an $\omega$-indexed inverse system of maps of 
level $\geq \alpha$ are also of level $\geq \alpha$;
and $(iii)$ if {\bf K} is
a co-elementary class, $k < \omega$ and 
$\mbox{Lev}_{\geq k}(\mbox{\bf K,K})=$ 
$\mbox{Lev}_{\geq k+1}(\mbox{\bf K,K})$, then 
$\mbox{Lev}_{\geq k}(\mbox{\bf K,K})=$ 
$\mbox{Lev}_{\geq\omega}(\mbox{\bf K,K})$.   

A space $X \in \mbox{\bf K}$ is {\it co-existentially closed in\/} {\bf K} if
$\mbox{Lev}_{\geq 0}(\mbox{\bf K},X)=$ 
$\mbox{Lev}_{\geq 1}(\mbox{\bf K},X)$.  We showed in \cite{Ban9} that
every infinite member of a co-inductive co-elementary class (such as
{\bf CH} itself, {\bf BS}, or the class {\bf CON} of continua) is a continuous
image of a space of the same weight that is co-existentially closed in
that class.  We show here that every compactum that is co-existentially   
closed in {\bf CON} (a {\it co-existentially closed continuum\/}) is both
indecomposable and of covering dimension one.
  
\end{abstract}

\noindent
{\it A.M.S. Subject Classification\/} (1991):
03C20, 54B35, 54C10, 54D05, 54D30, 54D80, 54F15, 54F45, 54F55.\\[10pt]


\section{Introduction.}\label{1}

This paper, a continuation of \cite{Ban2}--\cite{Ban9}, aims to carry
on the project of establishing model-theoretic concepts and methods
within the topological context; namely that of compacta (i.e., compact
Hausdorff spaces).  Since there is a precise duality between the
categories of compacta (plus continuous maps) and
commutative $B^*$-algebras (plus nonexpansive linear maps) 
(the Gel'fand-Na\u{\i}mark theorem \cite{Sim}),
our enterprise may also be seen as part of Banach model theory 
(see \cite{HH}--\cite{HI}).  The main difference is that we are doing
Banach model theory ``in the mirror,'' so to speak, and it is often the
case that a mirror can help one focus on features that might otherwise go
unnoticed.   

In the interests of being as self-contained as possible, we present
a quick review of our main tool, the topological ultracoproduct construction.
It is this construction, plus the landmark ultrapower theorem of
Keisler-Shelah \cite{CK}, that gets our project off the ground.
(Detailed accounts may be found in \cite{Ban2}--\cite{Ban8} and \cite{Gur}.)

We let {\bf CH} denote the category of compacta
and continuous maps.  In model theory, it is well known
that ultraproducts 
may be described in the language of category theory; i.e., as direct limits 
of (cartesian) products, where the directed 
set is the ultrafilter with reverse inclusion, and the system of products  
consists of cartesian products taken over the various sets in the ultrafilter.
(Bonding maps are just the obvious restriction maps.)  When we transport this
framework to the category-opposite of {\bf CH}, the result is the
{\bf topological ultracoproduct} (i.e., take an inverse limit of coproducts), 
and may be concretely
described as follows:  Given compacta $\langle X_i:i  \in I \rangle$ and
an ultrafilter $\cal D$ on $I$,
let $Y$ be the disjoint union $\bigcup_{i \in I}(X_i \times \{i\})$
(a locally compact space).  With $q:Y \to I$ the natural projection onto
the second co\"{o}rdinate (where $I$ has the discrete topology), we then
have the Stone-\v{C}ech lifting $q^{\beta}:\beta(Y)\to \beta(I)$.  Now the
ultrafilter $\cal D$ may be naturally viewed as an element of $\beta(I)$, 
and it is not hard to show that the topological ultracoproduct 
$\sum_{\cal D}X_i$ is the pre-image 
$(q^{\beta})^{-1}[\cal D]$.  
(The reader may be familiar with the Banach ultraproduct \cite{DK}. 
This construction is indeed the ultraproduct
in the category of Banach spaces and nonexpansive linear maps, and may
be telegraphically described using the recipe: take the usual ultraproduct,
throw away the infinite elements, and mod out by the subspace of 
infinitesimals.  Letting $C(X)$ denote the Banach space of continuous
real-valued (or complex-valued) continuous functions with $X$ as domain,
the Banach ultraproduct of $\langle C(X_i): i\in I\rangle$ via $\cal D$ is
just $C(\sum_{\cal D}X_i)$.)

If $X_i = X$
for all $i\in I$, then we have the {\bf topological ultracopower}
$XI\backslash {\cal D}$, a
subspace of $\beta(X\times I)$. In this case there is the Stone-\v{C}ech
lifting $p^{\beta}$ of the natural first-co\"{o}rdinate map 
$p:X\times I\to X$.  Its restriction to the ultracopower is a continuous
surjection, called the {\bf codiagonal map}, and is officially 
denoted $p_{X,{\cal D}}$ (with the occasional notation-shortening alias
possible). 
This map is dual to the natural diagonal map from a relational structure
to an ultrapower of that structure, and is 
not unlike the standard part map from nonstandard
analysis.)

When attention is restricted to the full subcategory {\bf BS} of Boolean 
spaces,  Stone duality assures us that the ultracoproduct construction
matches perfectly with the usual ultraproduct construction for Boolean 
lattices.  This says that ``dualized model theory'' in {\bf BS} is
largely a predictable rephrasing of the usual model theory of the elementary
class of Boolean
lattices.  In the category {\bf CH}, however, there is no similar match
(see \cite{Bana,Ros}); one is
forced to look for other (less direct) model-theoretic aids.  Fortunately  
there is a finitely axiomatizable AE Horn class of bounded distributive 
lattices,
the so-called {\bf normal disjunctive\/} lattices \cite{Ban8} (also called
Wallman lattices in \cite{Ban5}), comprising precisely the (isomorphic
copies of) {\bf lattice bases}, those lattices that serve as bases for the 
closed sets of compacta.  (To be more specific:  The normal disjunctive
lattices are precisely those bounded lattices $A$ such that there exists
a compactum $X$ and a meet-dense sublattice ${\cal A}$ of the closed set
lattice $F(X)$ of $X$ such that $A$ is isomorphic to $\cal A$.)  
We go from lattices to spaces, as in the case of Stone
duality, via the {\bf maximal spectrum\/} $S(\;)$, pioneered by H. Wallman 
\cite{Walm}.
$S(A)$ is the space of 
maximal proper filters of $A$; a typical basic closed set in $S(A)$ is the
set $a^{\sharp}$ of elements of $S(A)$ containing a given element $a \in A$.  
$S(A)$ is generally compact with this topology.
Normality, the condition that if
$a$ and $b$ are disjoint ($a \sqcap b = \bot$), then there are $a'$,
$b'$ such that $a \sqcap a' = b \sqcap b' = \bot$ and $a' \sqcup b' = \top$,
ensures that the maximal spectrum topology is Hausdorff.  Disjunctivity,
which says that for any two distinct lattice elements there is a nonbottom 
element that is below one
and disjoint from the other,
ensures that the map $a \mapsto a^{\sharp}$ takes
$A$ isomorphically onto the canonical closed set base for $S(A)$.  $S(\;)$
is contravariantly functorial: If $f:A \to B$ is a homomorphism of normal
disjunctive lattices and $M \in S(B)$, then $f^S(M)$ is the unique maximal
filter extending the prime filter $f^{-1}[M]$.  (For normal 
lattices, each prime filter is contained in a unique maximal one.)

The ultrapower theorem states that two relational structures are 
elementarily equivalent if and only if some ultrapower of one is isomorphic
to some ultrapower of the other.  One may easily extend this result, by
the use of added constant symbols, to show
that a function $f:A\to B$ between structures is an elementary
embedding if and only if there is an isomorphism of ultrapowers
$h:A^I/{\cal D} \to B^J/{\cal E}$ such that the obvious mapping square
commutes; i.e., such that $d_{\cal E}\circ f = h\circ d_{\cal D}$, where 
$d_{\cal D}$
and $d_{\cal E}$ are the natural diagonal embeddings.  This characterization
is used, in a thoroughly straightforward way, 
to {\it define\/} what it means for two compacta to be
{\bf co-elementarily equivalent} and for a map between compacta to be
a {\bf co-elementary map}.  It is a relatively easy
task to show, then, that $S(\;)$ converts ultraproducts
to ultracoproducts, elementarily equivalent lattices to co-elementarily
equivalent compacta, and elementary embeddings to 
co-elementary 
maps.  Furthermore, if $f:A\to B$ is a {\bf separative} embedding; i.e.,
an embedding such that
if $b \sqcap c = \bot$ in $B$, then there exists $a \in A$ such that
$f(a) \geq b$ and $f(a) \sqcap c = \bot$, then $f^S$ is a homeomorphism  
(see \cite{Ban2,Ban4,Ban5,Ban8,Gur}).  Because of this, there is much
flexibility in how we may obtain $\sum_{\cal D}X_i$:  Simply choose a
lattice base ${\cal A}_i$ for each $X_i$ and apply $S(\;)$ to the
ultraproduct $\prod_{\cal D}{\cal A}_i$.       

The spectrum functor falls far short of being a duality, except when restricted
to the Boolean lattices.  For this reason, one must take care not to jump
to too many optimistic conclusions; such as inferring that if compacta $X$ and 
$Y$
are co-elementarily equivalent, then there must
be lattice bases $\cal A$ for $X$ and $\cal B$ for $Y$ such that 
$\cal A \equiv \cal B$.  Similarly, one may not assume that a co-elementary
map is of the form $f^S$ for some elementary embedding.  This ``representation
problem'' has yet to be solved.\\
 
\section{An Ordinal-indexed Hierarchy of Maps.}\label{2}

Recall the definition of quantifier rank for first-order formulas in
prenex normal form:  $\varphi$ is of {\bf rank} 0 if it is quantifier
free; for $k<\omega$, $\varphi$ is of {\bf rank} $k+1$ if $\varphi$
is of the form $\pi\psi$, where $\psi$ is a prenex formula of rank $k$,
and $\pi$ is a prefix of like quantifiers, of polarity opposite to
that of the leading quantifier of $\psi$ (if there is one). 
We use the notation $\varphi(x_1,\dots ,x_n)$ to mean that the free variables
occurring in $\varphi$ come from the set
$\{x_1,...,x_n\}$.  For $k < \omega$, a
function $f:A \to B$ between structures is a {\bf map of level}
$\geq k$ if for every formula $\varphi(x_1\dots ,x_n)$ of rank $k$ 
and every
$n$-tuple $\langle a_1,\dots ,a_n\rangle$ from $A$, $A \models 
\varphi[a_1,\dots ,a_n]$ (if and) only if 
$B \models \varphi[f(a_1),\dots ,f(a_n)]$. (The obvious substitution
convention is being followed here.)  Maps of level $\geq 0$ are just
the embeddings; maps of level $\geq 1$ are often called {\it existential\/}
embeddings. (So the image under an existential embedding of one 
structure into another is existentially closed in the larger structure.)    
Of course an embedding is elementary if and only if it is of level 
$\geq \omega$; i.e., of level $\geq k$ for all $k<\omega$.

The following result is well known (see \cite{HSim}), and forms the basis 
upon which we
can export the model-theoretic notion of {\it map of level $k$\/}
to the topological context.
A function $f:A \to B$ between relational structures is a map of level
$\geq k+1$ if and only if
there are functions $g:A \to C$, $h:B \to C$
such that $g$ is an elementary embedding, $h$ is a map of level $\geq k$,
and  $g=h\circ f$.
($C$ may be taken to be an ultrapower of $A$, with $g$ the natural diagonal.)

We then define the notion of {\it map of level $k$\/} in the compact
Hausdorff context by use of this characterization.  $f:X \to Y$ is of
{\bf level} $\geq 0$ if it is a continuous surjection; for $k < \omega$,
$f$ is of {\bf level} $\geq k+1$ if there are functions $g:Z \to Y$,
$h:Z \to X$ such that $g$ is a co-elementary map, $h$ is a map of level
$\geq k$, and $g = f\circ h$.  If $f$ happens to be of level $\geq k$ for every
$k<\omega$, there is no obvious reason to infer that $f$ is co-elementary.
It therefore makes sense to carry the hierarchy into the transfinite, 
taking intersections at the limit stages and mimicking the inductive stage
above otherwise.  Thus we may talk of maps of level $\geq \alpha$ for
$\alpha$ any ordinal.  Clearly co-elementary maps are of level $\geq \alpha$
for each $\alpha$, but indeed there is no obvious assurance that the
converse is true.  The main goal of this section is to show that being
of level $\geq \omega$
is in fact equivalent to being co-elementary.\\

\subsection{Remarks.}\label{2.1}  
$(i)$ Co-elementary
equivalence is known \cite{Ban2,Ban5,Gur} to preserve important properties
of topological spaces, such as being infinite, being a continuum (i.e., 
connected), being Boolean (i.e., totally
disconnected), having (Lebesgue covering) dimension $n$, and being a
decomposable continuum.  If $f:X \to Y$ is a co-elementary map in {\bf CH},
then of course $X$ and $Y$ are co-elementarily equivalent ($X \equiv Y$). 
Moreover, since $f$ is a continuous surjection (see 
\cite{Ban2}), additional information about $X$ is transferred to $Y$.  For
instance, continuous surjections in {\bf CH} cannot raise {\bf weight\/} (i.e., 
the smallest cardinality of a possible topological base, and for many
reasons the right cardinal invariant to replace cardinality in the dualized
model-theoretic setting), so metrizability
(i.e., being of countable weight in the compact Hausdorff context) is
preserved.  Also local connectedness is preserved, since continuous surjections
in {\bf CH} are quotient maps.  Neither of these properties is an invariant
of co-elementary equivalence alone.  

$(ii)$ A number of properties, not generally preserved by continuous surjections
between compacta, are known \cite{Ban9} to be preserved by co-existential 
(i.e., level $\geq 1$) maps.  Among these are:
being infinite, 
being disconnected, having dimension $\leq n$, and  
being an (hereditarily) indecomposable continuum.\\

As is shown in \cite{Ban2}, co-elementary equivalence is an equivalence
relation (the sticking point being transitivity, of course), and the
composition of co-elementary maps is again a co-elementary map.  Furthermore,
there is the following ``closure under terminal factors'' property:  If
$f$ and $g\circ f$ are co-elementary maps, then so is $g$.  What makes these
(and many other) results work is the following lemma, an application of
a strong form of Shelah's version of the ultrapower theorem (see \cite{Ban9,
She}).
The following is a slight rephrasing of Lemma 2.1 in \cite{Ban9}, and is 
proved the same way.\\

\subsection{Lemma.}\label{2.2}
Let
$\langle\langle X_{\delta},f_{\delta},Y_{\delta}\rangle : \delta \in \Delta
\rangle$ be a family of triples, 
where $f_{\delta}$ either indicates co-elementary equivalence between
$X_{\delta}$ and $Y_{\delta}$, or is a co-elementary map from $X_{\delta}$
to $Y_{\delta}$, both spaces being compacta.
Then there is a
single ultrafilter witness to the fact. More precisely, there is an ultrafilter 
$\cal D$ on a
set $I$ and a family of homeomorphisms 
$\langle h_{\delta}: X_{\delta}I\backslash {\cal D} \to 
Y_{\delta}I\backslash {\cal D}: \delta \in \Delta \rangle$ 
such that $f_{\delta}\circ p_{X_{\delta},{\cal D}} = 
p_{Y_{\delta},{\cal D}}\circ h_{\delta}$ whenever $f_{\delta}$ is a 
co-elementary map.\\

In order for us to prove any substantial results concerning maps of level
$\geq \alpha$, we must extend the ultracoproduct construction from
compacta to continuous maps between compacta.  This was originally done in
\cite{Ban2}, but we need to establish some new facts about this construction.

Recall that if $f_i:X_i \to Y_i$ is a continuous map for each $i \in I$,
and $\cal D$ is an ultrafilter on $I$, then $\sum_{\cal D}f_i:
\sum_{\cal D}X_i \to \sum_{\cal D}Y_i$ may be defined as
$(\prod_{\cal D}f_i^F)^S$, where $f_i^F$ is just ``pulling closed
sets back to closed sets,'' and the ultraproduct map at the lattice level
is defined in the usual way.  When all the maps $f_i$ are equal to a
single map $f$, we have the ultracopower map, which we denote
$fI\backslash {\cal D}$.  
It is straightforward to show that 
ultracoproducts of continuous surjections (resp., homeomorphisms) are again
continuous surjections (resp., homeomorphisms).  In particular the
ultracopower operation $(\;)I\backslash {\cal D}$ is an endofunctor on
the category {\bf CH}.\\ 

\subsection{Proposition.}\label{2.3} 
Ultracoproducts of co-elementary maps are co-elementary maps.  More
specifically, 
if $\{i \in I: f_i\;\mbox{is a co-elementary map}\}\in {\cal D}$, then
$\sum_{\cal D}f_i$ is a co-elementary map also.\\

\noindent
{\bf Proof.} Let $f_i:X_i \to Y_i$ be given, $i \in I$, and set $J :=
\{i \in I: f_i\;\mbox{is a co-elementary map}\}\in {\cal D}$. 
We first consider the special case where $X_i := Y_iK_i\backslash {\cal E}_i$,
and $f_i := p_{{\cal E}_i}$ (the codiagonal map), for $i \in J$.  
Then $f_i = 
d_{{\cal E}_i}^S$, the image under the maximal spectrum functor of the
canonical diagonal embedding taking $F(Y_i)$ to $F(Y_i)^{K_i}\slash {\cal E}_i$.
Now each $d_{{\cal E}_i}$ is an elementary embedding; and an easy
consequence of the \L o\'{s} ultraproduct theorem is that ultraproducts
of elementary embeddings are elementary.  Since $S(\;)$ converts elementary
embeddings to co-elementary maps, we conclude that $\sum_{\cal D}f_i$ is
a co-elementary map in this case.  

In general, we have, for $i \in J$, homeomorphisms between ultracopowers
$h_i:X_iK_i\backslash {\cal D}_i \to Y_iL_i\backslash {\cal E}_i$, with
$f_i\circ p_{{\cal D}_i}= p_{{\cal E}_i}\circ h_i$.  When we take the 
ultracoproduct,
commutativity is preserved, $\sum_{\cal D}h_i$ is a homeomorphism, and
$\sum_{\cal D}p_{{\cal D}_i}$ and 
$\sum_{\cal D}p_{{\cal E}_i}$ are both co-elementary.  Thus 
$\sum_{\cal D}f_i$ is co-elementary, by closure under terminal factors.
$\dashv$\\  

The following analogue of \ref{2.3} can now be easily proved.\\

\subsection{Corollary.}\label{2.4} 
For each ordinal $\alpha$, ultracoproducts of maps of level $\geq \alpha$ 
are maps of level $\geq \alpha$.\\ 

\noindent
{\bf Proof.} The proof is by induction on $\alpha$.  Ultracoproducts of
continuous surjections are continuous surjections, so the result is 
established for $\alpha = 0$.  The inductive step at limit ordinals is
trivial, so it remains to prove the inductive step at successor ordinals.
But this follows immediately from \ref{2.3} and the definition of being of 
level $\geq \alpha + 1$. $\dashv$\\

Next we need closure under composition.\\

\subsection{Proposition.}\label{2.5}
For each ordinal $\alpha$, the composition of two maps of level $\geq \alpha$
is a map of level $\geq \alpha$.\\

\noindent
{\bf Proof.} Again we prove by induction on $\alpha$.  There is no problem
for $\alpha$ either zero or a positive limit ordinal, so assume the
composition of two maps of level $\geq \alpha$ is also of level $\geq \alpha$,
and let $f:X \to Y$ and $g:Y \to Z$ be maps of level $\geq \alpha + 1$.
By definition, there are maps $u:W \to X$, $v:W \to Y$ such that $u$ is
co-elementary, $v$ is of level $\geq \alpha$, and $u=f\circ v$.  By the
co-elementarity of $u$, there are ultracopowers 
$p:YI\backslash {\cal D}\to Y$,   
$q:WJ\backslash {\cal E}\to W$, and a homeomorphism    
$h:WJ\backslash {\cal E}\to 
YI\backslash {\cal D}$ such that $u\circ q = p\circ h$.  By our inductive 
hypothesis,
$v\circ q\circ h^{-1}$ is of level $\geq \alpha$; so we are justified in 
assuming
that the co-elementary part of a witness to a map's being of level $\geq
\alpha$ may be taken to be an ultracopower codiagonal map.

Getting back to $f$ and $g$, and using \ref{2.2}, there are maps
$p:YI\backslash {\cal D} \to Y$, 
$h:YI\backslash {\cal D} \to X$, 
$q:ZI\backslash {\cal D} \to Z$, 
$j:ZI\backslash {\cal D} \to Y$ such that $p$ and $q$ are codiagonal
maps, $h$ and $j$ are maps of level $\geq \alpha$, and the equalities
$p = f\circ h$ and $q = g\circ j$ both hold.  By \ref{2.4}, the ultracopower 
map $gI\backslash {\cal D}$ is of level $\geq \alpha + 1$.  Thus we have
further witnesses $u:W \to ZI\backslash {\cal D}$ and $v:W \to
YI\backslash {\cal D}$ such that $u$ is co-elementary, $v$ is of level 
$\geq \alpha$, and $u = (gI\backslash {\cal D})\circ v$.  Now $h\circ v$ is 
of
level $\geq \alpha$ by our inductive hypothesis, $q\circ v$ is co-elementary
by the long-established fact \cite{Ban2} that co-elementarity is closed
under composition, and it is a routine exercise to establish that  
$q\circ u = g\circ f\circ h\circ v$.  Thus $g\circ f$ is of level 
$\geq \alpha + 1$. $\dashv$\\

\subsection{Corollary.}\label{2.5.5}
Let $\alpha$ be an ordinal, $f:X\to Y$ a map of level $\geq \alpha+1$.  Then
there is an ultracopower map $p:YI\backslash {\cal D} \to Y$ and a map
$h:YI\backslash {\cal D}\to X$ of level $\geq \alpha$ such that 
$p = f\circ h$.\\

\noindent
{\bf Proof.} This is immediate from \ref{2.5}, plus the definitions of 
{\it co-elementary
map\/} and {\it map of level $\alpha+1$\/}. $\dashv$\\

\ref{2.5} gives us the following analogue of closure under terminal
factors for co-elementary maps.\\

\subsection{Corollary.}\label{2.6}
Let $\alpha$ be an ordinal.
If $h$ is of level $\geq \alpha$ and $f\circ h$ is of level $\geq \alpha + 1$,
then $f$ is of level $\geq \alpha + 1$.\\

\noindent
{\bf Proof.} Let $f:X \to Y$, $h:Z \to X$ be given, where $h$ is of level
$\geq \alpha$ and $g:=f\circ h$ is of level $\geq \alpha + 1$.  Then we have,
as witness to the level of $g$, maps $u:W \to Y$ and $v:W \to Z$ such 
that $u$ is co-elementary, $v$ is of level $\geq \alpha$, and $u=g\circ v$.
By \ref{2.5}, $h\circ v$ is of level $\geq \alpha$, and we have a witness to
the fact that $f$ is of level $\geq \alpha +1$. $\dashv$\\
 
We may now establish a needed consequence of \ref{2.3}, \ref{2.4},
and \ref{2.5}.\\

\subsection{Proposition.}\label{2.7}
Let $f:X \to Y$ be a map of level $\geq \omega$ between compacta. Then
there are maps $g:Z \to Y$ and $h:Z \to X$ such that $g$ is co-elementary,
$h$ is of level $\geq \omega$, and $g = f\circ h$.  Moreover, $g$ may be
taken to be an ultracopower codiagonal map.\\

\noindent
{\bf Proof.} For each $k<\omega$, $f$ is of level $\geq k+1$.  So let
$g_k:Z_k \to Y$ and $h_k:Z_k \to X$ witness the fact; each $g_k$ is
co-elementary, each $h_k$ is of level $\geq k$, and $g_k = f\circ h_k$.
Let $\cal D$ be any nonprincipal ultrafilter on $\omega$.  For each
$k<\omega$ we then have $\{k \in \omega: g_k\;\mbox{is co-elementary and}\;
h_k\;\mbox{is of level}\;\geq k\} \in {\cal D}$.  By \ref{2.3}, 
$\sum_{\cal D}g_k$ is co-elementary; by \ref{2.4}, $\sum_{\cal D}h_k$
is of level $\geq \omega$, and $\sum_{\cal D}g_k =
(f\omega\backslash {\cal D})\circ (\sum_{\cal D}h_k)$.  Let 
$p:X\omega\backslash {\cal D}
\to X$ and $q:Y\omega\backslash {\cal D}$ be the codiagonal maps.  Then, 
by \ref{2.5}, $p\circ (\sum_{\cal D}h_k)$ is of level $\geq \omega$.  We also
have the co-elementarity of $q\circ (\sum_{\cal D}g_k)$, as well as the equality
$f\circ p\circ (\sum_{\cal D}h_k) = q\circ (\sum_{\cal D}g_k)$; hence the 
desired result.  Further application of \ref{2.5} makes it possible to
arrange for $g:Z \to Y$ to be an ultracopower codiagonal map. $\dashv$\\

In order to prove that maps of level $\geq \omega$ are co-elementary, we
need a result on {\it co-elementary chains}.  Suppose  
$\langle X_n \stackrel{f_n}{\leftarrow} X_{n+1}: n < \omega \rangle$ is
an $\omega$-indexed inverse system of maps between compacta.
Then there is a compactum $X$ and
maps $g_n:X \to X_n$, $n < \omega$, such that the equalities 
$g_n = f_n\circ g_{n+1}$ all hold.  Moreover, $X$ is ``universal'' in the
sense that if $h_n:Y \to X_n$ is any other family of maps such that the
equations $h_n = f_n\circ h_{n+1}$ all hold, then there is a unique
$f:Y \to X$ such that $h_n = g_n\circ f$ for all $n<\omega$.  $X$ is the
{\bf inverse limit} of the sequence, and may be described as the
subspace $\{\langle x_0,x_1,\dots\rangle \in \prod_{n<\omega}X_n:
x_n = f_n(x_{n+1})\;\mbox{for each}\;n<\omega\}$.  The limit map $g_n$ is
then just the projection onto the $n$th factor.  

The inverse system is called
a {\bf co-elementary chain} if each $f_n$ is a co-elementary map.  We 
would like to conclude that, with co-elementary chains, the limit 
maps $g_n$ are also co-elementary.  This would give us a perfect analogue
of the Tarski-Vaught elementary chains theorem (see \cite{CK}).  In the
model-theoretic version, the proof uses induction on the complexity of
formulas, and is elegantly simple.  In our setting, however, it is not
entirely obvious how to proceed with a proof.  The result is still true,
but there is no simple elegant proof that we know of.  One proof is 
outlined in \cite{Ban9} (see Theorem 4.2 there).  It uses an elementary
chains analogue in Banach model theory, plus the Gel'fand-Na\u{\i}mark
duality theorem.  We present two more proofs in the next section; ones
that use only the techniques we have developed so far.  

An important step on the way to the co-elementary chains theorem is the
result that every map of level $\geq \omega$ is co-elementary.  It turns
out that this step itself uses the co-elementary chains theorem, but only
in a weak form.  Given a co-elementary chain
$\langle X_n \stackrel{f_n}{\leftarrow} X_{n+1}: n < \omega \rangle$,
we say the chain is {\bf representable} if there is an elementary chain
$\langle A_n \stackrel{r_n}{\rightarrow}  A_{n+1}: 
n < \omega \rangle$ of normal disjunctive lattices such that, 
for each $n < \omega$, $X_n = S(A_n)$ and $f_n = r_n^S$.\\  
 
\subsection{Lemma.}\label{2.8}
Let $\langle X_n \stackrel{f_n}{\leftarrow} X_{n+1}: n < \omega \rangle$ 
be a representable co-elementary chain, with inverse limit $X$ and
limit maps $g_n:X \to X_n$, $n<\omega$.  Then Each $g_n$ is a 
co-elementary map.\\

\noindent
{\bf Proof.}  Let  
$\langle A_n \stackrel{r_n}{\rightarrow}  A_{n+1}: 
n < \omega \rangle$ represent our co-elementary chain in the sense given
above.  Let $A$ be the direct limit of this direct system, with limit 
maps $t_n:A_n \to A$.  Then the Tarski-Vaught theorem says that each
$t_n$ is an elementary embedding.  Since the maximal spectrum functor
converts elementary embeddings to co-elementary maps, we then have 
$Y:= S(A)$, and co-elementary maps $h_n := t_n^S$.  Let $f:Y \to X$
be defined by the equalities $h_n = g_n\circ f$.  Applying the closed
set functor $F(\;)$ to the representing elementary chain, letting
$u_n:A_n \to F(X_n)$ be the natural separative embedding, we have
the embeddings $g_n^F\circ u_n: A_n \to F(X)$.  We then get $u:A \to
F(X)$, defined by $u\circ t_n = g_n^F\circ u_n$.  Let $g:= u^S:X \to Y$.
Then we have, applying $S(\;)$, and noting that the maps $u_n^S$ are
canonical homeomorphisms, $u_n^S\circ g_n = h_n\circ g$.  This implies that
$f$ and $g$ are inverses of one another; hence that the maps $g_n$ are
co-elementary. $\dashv$\\   

We are now ready to prove the main result of this section.\\

\subsection{Theorem.}\label{2.9}
Every map of level $\geq \omega$ is co-elementary.\\

\noindent
{\bf Proof.}  Let $f_0:X_0 \to Y_0$ be a map of level $\geq \omega$.  We 
build a ``co-elementary ladder'' over this map as follows: 
By \ref{2.7}, there are maps $g_0:Y_1 \to Y_0$ and $h_0:Y_1 \to X_0$ such
that $g_0$ is co-elementary, $h_0$ is of level $\geq \omega$, and
$g_0 = f_0\circ h_0$.  Moreover, we may (and do) take $g_0$ to be an
ultracopower codiagonal map.  
Since $h_0$ is of level $\geq \omega$, we have maps
$j_0:X_1 \to X_0$ and $f_1:X_1 \to Y_1$ such that $j_0$ is co-elementary,
$f_1$ is of level $\geq \omega$, and $j_0 = h_0\circ f_1$.  As before,
we take $j_0$ to be an ultracopower codiagonal map.  
This completes
the first ``rung'' of the ladder, and we repeat the process for the
map $f_1:X_1 \to Y_1$.  In the end, we have two co-elementary chains
$\langle X_n \stackrel{j_n}{\leftarrow} X_{n+1}: n < \omega \rangle$ and 
$\langle Y_n \stackrel{g_n}{\leftarrow} Y_{n+1}: n < \omega \rangle$, 
with inverse limits $X$ and $Y$ respectively.  For each $n < \omega$, let
$v_n:X \to X_n$ and $w_n:Y \to Y_n$ be the limit maps, defined by the
equalities $v_n = j_n\circ v_{n+1}$, $w_n = g_n\circ w_{n+1}$.   

Now each successive entry is an ultracopower of the last; hence these  
co-elementary chains are representable (by elementary chains of iterated
ultrapowers).   By \ref{2.8}, the maps $v_n$ and $w_n$ are co-elementary.
  
Consider now the maps
$h_n:Y_{n+1}\to X_n$, $n<\omega$.  These, along with the maps $f_n$, give
rise to the existence of maps $f:X \to Y$
and $h:Y \to X$ that are unique with the property that for all $n < \omega$,
$w_n\circ f = f_n\circ v_n$ and $v_n\circ h = h_n\circ w_{n+1}$.  The
uniqueness feature ensures that $f$ and $h$ are inverses of one another;
thus $f_0$ is co-elementary, by closure under terminal factors. $\dashv$\\
 
\section{Inverse Limits of $\alpha$-chains.}\label{3}

In this section we prove the co-elementary chains theorem in two different
ways, both of which use \ref{2.9}.

If $\alpha$ is an ordinal, an inverse system  
$\langle X_n \stackrel{f_n}{\leftarrow} X_{n+1}: n < \omega \rangle$ 
of maps between compacta
is an $\alpha$-{\bf chain} if each $f_n$ is a map of level $\geq \alpha$. 
By the $\alpha$-{\bf chains theorem}, we mean the statement that the limit
maps of every $\alpha$-chain are maps of level $\geq \alpha$.  
(So, for example, the 0-chains theorem is a well-known exercise.)
Because of
\ref{2.9}, the co-elementary chains theorem is just the $\omega$-chains
theorem; and this case clearly follows from the conjunction of the cases
$\alpha < \omega$.  While we do ultimately prove the $\alpha$-chains theorem
for general $\alpha$, we first take
a slight detour and establish the $\alpha = \omega$ case separately.  The
main reason for doing this (aside from the fact that we discovered this
case first in an abortive attempt to establish the general case) is that it
uses the following result, which is of further use later on, as well as
being of some independent interest.\\    

\subsection{Lemma.}\label{3.1}
Let $f:X \to Y$ be a function between compacta, let $\alpha$
be an ordinal, and let ${\cal B}$ be a lattice base for $Y$. 
Suppose that for each finite $\delta \subseteq {\cal B}$ there
is a map $g_{\delta}:Y \to Z_{\delta}$, of level $\geq \alpha$, such that
$g_{\delta}\circ f$ is of level $\geq \alpha$ and for each $B \in \delta$,
$g_{\delta}^{-1}[g_{\delta}[B]]=B$ (i.e., $B$ is 
$g_{\delta}$-{\it saturated\/}).  
Then $f$ is a map of level $\geq \alpha$.\\

\noindent
{\bf Proof.}  The proof below uses the basic idea for proving Theorem 3.3
in \cite{Ban9}.  

For each ordinal $\alpha$, let $\mbox{\bf A}_{\alpha}$ be the assertion of
the lemma for maps of level $\geq \alpha$.  Then 
$\mbox{\bf A}_{\omega}$ follows immediately from
the conjunction of the assertions $\mbox{\bf A}_{\alpha}$ for $\alpha$
finite.  In view of \ref{2.9}, then, we may focus our attention on the
finite case.  While our proof is not by induction, it does require a separate
argument for the case $\alpha = 0$.

Let $B \in {\cal B}$.  If $\delta \supseteq \{B\}$, then $B$ is 
$g_{\delta}$-saturated; so $f^{-1}[B] = 
f^{-1}[g_{\delta}^{-1}[g_{\delta}[B]]]=
[g_{\delta}\circ f]^{-1}[g_{\delta}[B]]$, a closed subset of $X$.  Thus $f$
is continuous.  Suppose $f$ fails to be surjective.  Then, because $f$ is
continuous, we have
disjoint nonempty $B,C \in {\cal B}$ with $f[X]\subseteq B$.  Pick
$\delta \supseteq \{B,C\}$.  Then both $B$ and $C$ are $g_{\delta}$-saturated;
hence $g_{\delta}[B]$, and  
$g_{\delta}[C]$ are nonempty and disjoint.  But then $g_{\delta}\circ f$
fails to be surjective.  This establishes $\mbox{\bf A}_0$.  

In the sequel we fix $\alpha < \omega$, and prove the
assertion $\mbox{\bf A}_{\alpha +1}$.

Let $\Delta$ be the set of finite subsets of $\cal B$.  Using \ref{2.2}, there
is a single ultrafilter $\cal D$ on a set $I$ that may be used to 
witness the hypothesis of $\mbox{\bf A}_{\alpha +1}$.  To be precise, for each  
$\delta \in
\Delta$, the mapping diagram $\mbox{\bf D}_{\delta}$ consists of
maps $p_{\delta}$, $h_{\delta}$, $k_{\delta}$,
from $Z_{\delta}I\backslash {\cal D}$ to $Z_{\delta}$, $Y$, and $X$
respectively, such that $p_{\delta}$ is the codiagonal map (so co-elementary),
$h_{\delta}$ and $k_{\delta}$ are each of level $\geq \alpha$, and 
$p_{\delta} = g_{\delta}\circ h_{\delta} = g_{\delta}\circ f\circ k_{\delta}$.  
To this diagram we adjoin the codiagonal map $q:YI\backslash {\cal D}
\to Y$, and define $r:= k_{\delta}\circ (g_{\delta}I\backslash {\cal D}):
YI\backslash {\cal D} \to X$.  By \ref{2.4} and \ref{2.5}, $r$ is a map
of level $\geq \alpha$; we would be done, therefore, if the equality
$q = f\circ r$ were true.  Not surprisingly, this equality is generally false.
What {\it is\/} true are the equalities 
$g_{\delta}\circ q = g_{\delta} \circ f\circ r$.  To take advantage of this,   
we form an ``ultracoproduct'' of the diagrams 
$\mbox{\bf D}_{\delta}$.  

For each $\delta \in \Delta$, let $\hat{\delta}:= \{\gamma \in \Delta:
\delta \subseteq \gamma\}$.  Then the set $\{\hat{\delta}: \delta \in
\Delta\}$ clearly satisfies the finite intersection property, and hence
extends to an ultrafilter ${\cal H}$ on $\Delta$.  Form the 
``${\cal H}$-ultracoproduct'' diagram {\bf D} in the obvious way.
Then we have the codiagonal maps $u:X\Delta\backslash {\cal H} \to X$ and
$v:Y\Delta\backslash {\cal H} \to Y$.
Moreover, again by \ref{2.4} and \ref{2.5},
$u\circ (r\Delta\backslash {\cal H})$ is of level $\geq \alpha$.  We
will be done, therefore, once we show that 
$f\circ u\circ (r\Delta\backslash {\cal H})= 
v\circ (q\Delta\backslash {\cal H})$.

Now the map on the left is just $v\circ ((f\circ r)\Delta\backslash {\cal H})$.
Suppose $x \in (YI\backslash {\cal D})\Delta\backslash {\cal H}$ is sent to
$y_1$ under the left map and to $y_2$ under the right.  Let $y_1':=
[(f\circ r)\Delta\backslash {\cal H}](x)$ and
$y_2':=
[q\Delta\backslash {\cal H}](x)$.  Then $[\sum_{\cal H}g_{\delta}](y_1')= 
[\sum_{\cal H}g_{\delta}](y_2')$.  Assume $y_1 \neq y_2$.  Then, by the 
nature of codiagonal maps, there exist disjoint $B_1,B_2 \in {\cal B}$,
containing $y_1$ and $y_2$ in their respective interiors, such that
$B_1^{\Delta}/{\cal H} \in y_1'$ and 
$B_2^{\Delta}/{\cal H} \in y_2'$.  If $\delta \supseteq \{B_1,B_2\}$,
then both $B_1$ and $B_2$ are $g_{\delta}$-saturated.  Thus
$\{\delta \in \Delta: g_{\delta}[B_1] \cap g_{\delta}[B_2] = \emptyset\}
\in {\cal H}$; hence $\prod_{\cal H}g_{\delta}[B_1]$ and
$\prod_{\cal H}g_{\delta}[B_2]$ are disjoint subsets of 
$\prod_{\cal H}F(Z_{\delta})$.   
Now $\prod_{\cal H}g_{\delta}[B_1] \in [\sum_{\cal H}g_{\delta}](y_1')$
and $\prod_{\cal H}g_{\delta}[B_2] \in [\sum_{\cal H}g_{\delta}](y_2')$;
from which we conclude that
$[\sum_{\cal H}g_{\delta}](y_1') \neq
[\sum_{\cal H}g_{\delta}](y_2')$.  This contradiction tells us that $y_1
= y_2$ after all, completing the proof. $\dashv$\\

We can now give a new proof of the co-elementary chains theorem
(Theorem 4.2 in \cite{Ban9}), one where no Banach model theory is used.

\subsection{Theorem.}\label{3.2} 
Let $\langle X_n \stackrel{f_n}{\leftarrow} X_{n+1}: n < \omega \rangle$  
be a co-elementary chain of compacta, with inverse limit $X$.
Then the limit maps $g_n:X \to X_n$, $n<\omega$, are all co-elementary.\\

\noindent
{\bf Proof.}
We first prove a weak version of the theorem.  This version appears as
Proposition 4.1 in \cite{Ban9}.  
Let $\langle X_n \stackrel{f_n}{\leftarrow} X_{n+1}: n < \omega \rangle$  
be a co-elementary chain of compacta.  Then there exists a compactum $Y$
and co-elementary maps $h_n:Y \to X_n$, $n < \omega$, such that all the
equalities $h_n = f_n\circ h_{n+1}$ hold.  The proof of this is quite
easy, and we repeat it here for the sake of completeness.

By \ref{2.2}, there is an ultrafilter ${\cal D}$ on a set $I$ and 
homeomorphisms $k_n: X_{n+1}I\backslash {\cal D} \to X_nI\backslash {\cal D}$,
$n<\omega$, such that all the equalities $p_n\circ k_n = f_n\circ p_{n+1}$
hold (where the maps $p_n$ are the obvious codiagonals).  Let $Y$ be
the inverse limit of this system, with limit maps $j_n:Y \to
X_nI\backslash {\cal D}$.  Since each $k_n$ is a homeomorphism, so is each
$j_n$, and we set $h_n:= p_n\circ j_n$, a co-elementary map.  Clearly
$f_n\circ h_{n+1} = h_n$ always holds, and there is a map $h:Y \to X$,
uniquely defined by the equalities $g_n\circ h = h_n$.

Now consider the chain of embeddings  
$\langle F(X_n) \stackrel{f_n^F}{\rightarrow} F(X_{n+1}): n < \omega \rangle$,
with direct limit $\cal A$, and limit embeddings $r_n: F(X_n) \to {\cal A}$.
Then (see the argument in \ref{2.8}) we may treat $X$ as $S({\cal A})$
and each $g_n$ as $r_n^S$.  (Note: we cannot hope for these embeddings to
be elementary.)  For each finite $\delta \subseteq {\cal A}$, there is 
a least $n_{\delta} < \omega$ such that each member of $\delta$ is in the
range of $r_n$ for $n\geq n_{\delta}$.  This tells us that $X$ has a
lattice base $\cal A$ such that for each finite $\delta \subseteq {\cal A}$
and each $A \in \delta$, $A$ is $g_{n_{\delta}}$-saturated.  This puts is
in a position to use \ref{3.1}.

We prove that each $g_n$ is of level $\geq \alpha$, for $\alpha < \omega$,
by induction on $\alpha$.  Clearly each $g_n$ is of level $\geq 0$; so 
assume each $g_n$ to be of fixed level $\geq \alpha$.  Then, by \ref{3.1},
$h$ is of level $\geq \alpha$ too.  Since each $h_n$ is co-elementary,
we have now a witness to the fact that each $g_n$ is of level 
$\geq \alpha +1$.  Thus each $g_n$ is of level $\geq \omega$, and is 
hence co-elementary by \ref{2.9}.  $\dashv$\\   

We had originally thought that \ref{3.1} could be used to prove the
$\alpha$-chains theorem in general,
but were unable to get our idea to work.  What is missing is a weak
version of the assertion, namely the existence of a compactum $Y$ and
maps $h_n:Y \to X_n$ of level $\geq \alpha$ such that all the equalities
$h_n= f_n\circ h_{n+1}$ hold.  If we could do this, then we could prove
the strong version by induction on finite $\alpha$:  The $\alpha = 0$ case
is known; assuming the assertion true for fixed $\alpha$, and
that we are given an $(\alpha +1)$-chain, we find our compactum $Y$ and
maps $h_n$, all of level $\geq \alpha +1$.  The maps $g_n$ are of level
$\geq \alpha$ by the inductive hypothesis, and we conclude that $h$ is
of level $\geq \alpha$, by \ref{3.1}.   Then each $g_n$ is of level
$\geq \alpha +1$, by \ref{2.6}.

Rather than pursue the tack just outlined, we abandon \ref{3.1} in favor
of a similar-sounding (but somewhat different) lemma.\\

\subsection{Lemma.}\label{3.3}
Let $f:X \to Y$ be a function between compacta, let $\alpha$
be an ordinal, and let ${\cal A}$ be a lattice base for $X$. 
Suppose that for each finite $\delta \subseteq {\cal A}$ there
is a map $g_{\delta}:X \to Z_{\delta}$, of level $\geq \alpha$, and a map
$h_{\delta}:Z_{\delta} \to Y$, of level $\geq \alpha +1$, such that
$f = h_{\delta}\circ g_{\delta}$, and each member of $\delta$ is 
$g_{\delta}$-saturated.
Then $f$ is a map of level $\geq \alpha +1$.\\

\noindent
{\bf Proof.}
Assume that $f:X \to Y$, $\cal A$, and $\alpha$ are
fixed, with $\Delta$ the set of all finite subsets of $\cal A$.  The
ultrafilter $\cal H$ on $\Delta$ is exactly as in \ref{3.1}.
For each $\delta \in \Delta$, the
diagram $\mbox{\bf D}_{\delta}$ consists of continuous surjections
$g_{\delta}:X \to Z_{\delta}$,
$h_{\delta}:Z_{\delta} \to Y$, a codiagonal map $p:YI\backslash {\cal D} \to
Y$, and a continuous surjection $k_{\delta}:YI\backslash {\cal D} \to 
Z_{\delta}$. ($\cal D$ need not depend on $\delta$, by \ref{2.2}, but that
fact is not essential to the argument.)  The maps $g_{\delta}$ and   
$k_{\delta}$ are of level $\geq \alpha$, and the
equalities $f = h_{\delta}\circ g_{\delta}$ and $p = h_{\delta}\circ
k_{\delta}$ both hold.

We form the ``ultracoproduct'' diagram as in \ref{3.1}, adding the codiagonal
maps $u:X\Delta\backslash {\cal H} \to X$,
$v:Y\Delta\backslash {\cal H} \to Y$, along with our original map $f$.
We then define the relation $j := u\circ (\sum_{\cal H}g_{\delta})^{-1}\circ 
(\sum_{\cal H}k_{\delta}): (YI\backslash {\cal D})\Delta\backslash {\cal H}
\to X$.  Once we show $j$ is a map of level $\geq \alpha$, and that
$f\circ j = v\circ (p\Delta \backslash {\cal H})$, we will have a witness
to the fact that $f$ is of level $\geq \alpha +1$.   

To show $j$ is a function, it suffices to show that the kernel of
$\sum_{\cal H}g_{\delta}$ is contained within the kernel of $u$.  Indeed,
suppose $x_1, x_2 \in X\Delta\backslash {\cal H}$ are such that 
$u(x_1) \neq u(x_2)$.  Then there are disjoint $A_1, A_2 \in {\cal A}$, 
containing $u(x_1)$ and $u(x_2)$ in their respective interiors, such that
$A_1^{\Delta}\slash {\cal H} \in x_1$ and
$A_2^{\Delta}\slash {\cal H} \in x_2$.  If $\delta \supseteq \{A_1,A_2\}$,
then both $A_1$ and $A_2$ are $g_{\delta}$-saturated, so
$\{\delta \in \Delta: g_{\delta}[A_1]\cap g_{\delta}[A_2] = \emptyset\}\in
{\cal H}$.  Thus $\prod_{\cal H}g_{\delta}[A_1]$ and 
$\prod_{\cal H}g_{\delta}[A_2]$ are disjoint subsets of 
$\prod_{\cal H}F(Z_{\delta})$, and are elements of 
$[\sum_{\cal H}g_{\delta}](x_1)$ and 
$[\sum_{\cal H}g_{\delta}](x_2)$, respectively.  Thus
$[\sum_{\cal H}g_{\delta}](x_1) \neq
[\sum_{\cal H}g_{\delta}](x_2)$, so $j$ is a function.  That $j$ is surjective
is clear; that $f\circ j = v\circ (p\Delta)\backslash {\cal H})$ is a
simple diagram chase.  Since $j^{-1} = (\sum_{\cal H}k_{\delta})^{-1}
\circ (\sum_{\cal H}g_{\delta})\circ u^{-1}$, and $\sum_{\cal H}g_{\delta}$
is a closed map, we conclude that $j$ is continuous.  Now  
$\sum_{\cal H}k_{\delta}$ and
$\sum_{\cal H}g_{\delta}$ are maps of level $\geq \alpha$, by \ref{2.4},
and $u\circ(\sum_{\cal H}g_{\delta})^{-1}$ is of level $\geq \alpha +1$,
by \ref{2.6}.  Thus $j$ is of level $\geq \alpha$, by \ref{2.5}. $\dashv$\\

We are now ready to establish the $\alpha$-chains theorem in general.\\

\subsection{Theorem.}\label{3.4} 
Let $\alpha$ be a fixed ordinal, and
let $\langle X_n \stackrel{f_n}{\leftarrow} X_{n+1}: n < \omega \rangle$  
be an $\alpha$-chain of compacta, with inverse limit $X$.  Then the limit
maps $g_n:X \to X_n$, $n<\omega$, are all of level $\geq \alpha$.\\

\noindent
{\bf Proof.}
Use induction on $\alpha$.  As mentioned above, we need only consider
finite $\alpha$, and the $\alpha = 0$ case is an easy exercise.  So 
assume the $\alpha$-chains theorem to be true for some fixed $\alpha$, and let
$\langle X_l \stackrel{f_l}{\leftarrow} X_{l+1}: l < \omega \rangle$  
be an $(\alpha + 1)$-chain.  Fix $n < \omega$.  With the aim of applying
\ref{3.3}, $Y$ is $X_n$, and $f$ is $g_n$.  As in the proof of \ref{3.2},
$\cal A$ is the direct limit of the system
$\langle F(X_l) \stackrel{f_l^F}{\rightarrow} F(X_{l+1}): l < \omega \rangle$
of normal disjunctive lattices.  Given finite $\delta \subseteq {\cal A}$,
there is some (least) $m > n$ such that each member of $\delta$ is 
$g_m$-saturated.  Let $Z_{\delta}$ and $g_{\delta}$ be $X_m$ and $g_m$,
respectively, with $h_{\delta}$ the obvious finite composition of the
maps $f_k$, as $k$ runs from $n$ to $m-1$.
$g_{\delta}$ is of level $\geq \alpha$ by our induction hypothesis;   
$h_{\delta}$ is of level $\geq \alpha + 1$ by \ref{2.5}.
By \ref{3.3}, then, $g_n$ is a map of level $\geq \alpha +1$. $\dashv$\\

\section{When Levels Collapse.}\label{4}

Here we address the issue of when there is a collapsing of levels of maps
between classes of compacta.  Let {\bf K} and {\bf L} be subclasses of
{\bf CH}, and define 
$\mbox{Lev}_{\geq \alpha}(\mbox{\bf K},\mbox{\bf L})$  to be the class of
maps of level $\geq \alpha$, with domains in {\bf K} and ranges in {\bf L}. 
(If one of the classes happens to be a single homeomorphism type, say
{\bf K} is the homeomorphism type of $X$, then we write 
$\mbox{Lev}_{\geq \alpha}(X,\mbox{\bf L})$ to simplify notation. (Etc.))
Recall that a class {\bf K} is a {\bf co-elementary class} if {\bf K} is
closed under ultracoproducts and co-elementary equivalence.  (Of course,
being closed under co-elementary equivalence is tantamount to being closed
under ultracopowers and co-elementary images; so we could replace the
criteria for being a co-elementary class with the conditions of being
closed under ultracoproducts and co-elementary images.)

The first result of this section is reminiscent of Robinson's test from
model theory, and its proof is very similar to that of \ref{2.9}.\\

\subsection{Theorem.}\label{4.1}
Suppose  {\bf K} and {\bf L} are closed under ultracopowers,  
that $\alpha < \omega$, and that
$\mbox{Lev}_{\geq \alpha}(\mbox{\bf K},\mbox{\bf L})=
\mbox{Lev}_{\geq \alpha +1}(\mbox{\bf K},\mbox{\bf L})$ and
$\mbox{Lev}_{\geq \alpha}(\mbox{\bf L},\mbox{\bf K})=
\mbox{Lev}_{\geq \alpha +1}(\mbox{\bf L},\mbox{\bf K})$.  Then
$\mbox{Lev}_{\geq \alpha}(\mbox{\bf K},\mbox{\bf L})=
\mbox{Lev}_{\geq \omega}(\mbox{\bf K},\mbox{\bf L})$ and
$\mbox{Lev}_{\geq \alpha}(\mbox{\bf L},\mbox{\bf K})=
\mbox{Lev}_{\geq \omega}(\mbox{\bf L},\mbox{\bf K})$. \\ 

\noindent
{\bf Proof.}
Let $f_0:X_0 \to Y_0$ be a map of level $\geq \alpha$ from a member of
{\bf K} to a member of {\bf L}.  Then we build a ``co-elementary ladder,''
similar to the one in the proof of \ref{2.9}, as follows:

Since $f_0$ is also of level $\geq \alpha +1$, 
there are maps $g_0:Y_1 \to Y_0$ and $h_0:Y_1 \to X_0$ such
that $g_0$ is co-elementary, $h_0$ is of level $\geq \alpha$, and
$g_0 = f_0\circ h_0$.  Moreover, we may (and do) take $g_0$ to be an
ultracopower codiagonal map; so, in particular, $Y_1 \in \mbox{\bf L}$, and  
$h_0$ is of level $\geq \alpha +1$.  Thus we have maps
$j_0:X_1 \to X_0$ and $f_1:X_1 \to Y_1$ such that $j_0$ is co-elementary,
$f_1$ is of level $\geq \alpha$, and $j_0 = h_0\circ f_1$.  As before,
we take $j_0$ to be an ultracopower codiagonal map, so $X_1 \in \mbox{\bf K}$.  
This completes
the first ``rung'' of the ladder, and we repeat the process for the
map $f_1:X_1 \to Y_1$, a map of level $\geq \alpha +1$.  

The rest of the proof proceeds exactly like the proof of \ref{2.9}, and
we conclude that $f_0$ is co-elementary.  $\dashv$\\

With the aid of \ref{3.1}, \ref{4.1} has some interesting variations.  We first
restate what in \cite{Ban9} we call the ``sharper'' L\"{o}wenheim-Skolem
theorem.  In the sequel, $w(X)$ stands for the weight of a space $X$.\\

\subsection{Theorem.}\label{4.2} (Theorem 3.1 of \cite{Ban9})  Let
$f:X \to Y$ be a continuous surjection between compacta, with $\kappa$
an infinite cardinal such that $w(Y) \leq \kappa \leq w(X)$.  Then there
is a compactum $Z$ and continuous surjections $g:X \to Z$, $h:Z \to Y$
such that $w(A) = \kappa$, $g$ is a co-elementary map, and $f = h\circ g$.\\

We next bring 3.1 into the picture with the following strengthening of
Theorem 3.3 in \cite{Ban9}.\\

\subsection{Theorem.}\label{4.3} Let $f:X \to Y$ be a function between
compacta, let $\alpha$ be an ordinal, and let $\kappa \leq w(Y)$ be an
infinite cardinal.  Suppose that for each compactum $Z$ of weight $\kappa$,
and each co-elementary map $g:Y \to Z$, the composition $g\circ f$ is a 
map of level $\geq \alpha$.  Then $f$ is a map of level $\geq \alpha$.\\

\noindent
{\bf Proof.} We let $\Delta$ be the set of finite subsets of $F(Y)$.
For each $\delta \subseteq \Delta$ there is a countable elementary
sublattice ${\cal A}_{\delta}$ of $F(Y)$, with $\delta \subseteq
{\cal A}_{\delta}$.  Let $W_{\delta}:= S({\cal A}_{\delta})$ (a space of
weight $\aleph_0$), with $r_{\delta}:Y \to W_{\delta}$ denoting the
co-elementary map that arises from the inclusion ${\cal A}_{\delta}
\subseteq F(Y)$.  Then every member of $\delta$ is $r_{\delta}$-saturated.
By \ref{4.2}, there is a compactum $Z_{\delta}$ of weight $\kappa$, and
continuous surjections $g_{\delta}: Y \to Z_{\delta}$, $t_{\delta}:
Z_{\delta} \to W_{\delta}$, such that $g_{\delta}$ is co-elementary and
$r_{\delta} = t_{\delta}\circ g_{\delta}$.  So each $g_{\delta}$ is a 
co-elementary map onto a compactum of weight $\kappa$; by hypothesis, then,
$g_{\delta}\circ f$ must be a map of level $\geq \alpha$.  By \ref{3.1},
$f$ must be a map of level $\geq \alpha$. $\dashv$\\

For any class {\bf K} and cardinal $\kappa$, let $\mbox{\bf K}_{\kappa}
:= \{X \in \mbox{\bf K}: w(X) = \kappa\}$.  The following is a variation
(though not, strictly speaking, an improvement) on \ref{4.1}.   \\

\subsection{Theorem.}\label{4.4}
Suppose  {\bf K} and {\bf L} are closed under ultracopowers, as well as
co-elementary images,  
that $0<\alpha < \omega$, and, for some infinite cardinal $\kappa$, that
$\mbox{Lev}_{\geq \alpha}(\mbox{\bf K}_{\kappa},\mbox{\bf L}_{\kappa})=
\mbox{Lev}_{\geq \alpha +1}(\mbox{\bf K}_{\kappa},\mbox{\bf L}_{\kappa})$ and
$\mbox{Lev}_{\geq \alpha}(\mbox{\bf L}_{\kappa},\mbox{\bf K}_{\kappa})=
\mbox{Lev}_{\geq \alpha +1}(\mbox{\bf L}_{\kappa},\mbox{\bf K}_{\kappa})$.  Then
$\mbox{Lev}_{\geq \alpha}(\mbox{\bf K},\mbox{\bf L})=
\mbox{Lev}_{\geq \omega}(\mbox{\bf K},\mbox{\bf L})$ and
$\mbox{Lev}_{\geq \alpha}(\mbox{\bf L},\mbox{\bf K})=
\mbox{Lev}_{\geq \omega}(\mbox{\bf L},\mbox{\bf K})$. (The assertion also
holds in the case $\alpha = 0$, if we assume that neither {\bf K} nor
{\bf L} contains any finite spaces.) \\ 

\noindent
{\bf Proof.} Let $f:X \to Y$ be a map of level $\geq \alpha$, between 
members of {\bf K} and {\bf L} respectively.  By \ref{4.1}, it suffices
to show that $f$ is of level $\geq \alpha +1$.  Assume first that 
$\kappa \leq w(Y)$.  By \ref{4.3}, it suffices to show that for each compactum
$Z$ of weight $\kappa$ and each co-elementary map $g:Y \to Z$, we have
that $g\circ f$ is of level $\geq \alpha +1$.  So let $g:Y \to Z$ be given.
By \ref{4.2}, there is a factorization $u: X \to W$, $v:W \to Z$ such 
that $u$ is co-elementary, $w(W) = w(Z) = \kappa$, and $g\circ f=
v\circ u$.  Now $W \in \mbox{\bf K}_{\kappa}$ and
$Z \in \mbox{\bf L}_{\kappa}$, and $g\circ f$ is of level $\geq \alpha$.
Thus, by \ref{2.6}, $v$ is also of level $\geq \alpha$.  By hypothesis,
$v$ is of level $\geq \alpha +1$; consequently, so is $g\circ f$. 

If $\kappa > w(Y)$, and we are dealing with the case $\alpha > 0$, then
we must consider the possibility that $Y$ is finite.  But $f$ is a 
co-existential map, and hence clearly a bijection (i.e., a homeomorphism)
in that situation.
So we may as well assume that $Y$ is infinite.  If we are dealing with
the case $\alpha = 0$, then we take $Y$ to be infinite by {\it fiat\/}.

That said, we find an ultrafilter $\cal D$ on a set $I$ such that
$w(YI\backslash {\cal D}) \geq \kappa$ (see \cite{Ban2}).  By the argument
in the first paragraph, since both {\bf K} and {\bf L} are closed under
ultracopowers, we conclude that $fI\backslash {\cal D}$ is of level
$\geq \alpha +1$.  From our work in \S2, we infer that $f$ is of level
$\geq \alpha +1$ too. $\dashv$.\\  

Given an ordinal $\alpha$, we say $X \in \mbox{\bf K}$ is 
$\alpha$-{\bf closed in K} if 
$\mbox{Lev}_{\geq 0}(\mbox{\bf K},X)= 
\mbox{Lev}_{\geq \alpha}(\mbox{\bf K},X)$.  (1-closed = co-existentially
closed \cite{Ban9}.)  Define $\mbox{\bf K}^{\alpha}:= \{X \in \mbox{\bf K}:
X\;\mbox{is}\;\alpha\mbox{-closed in \bf K}\}$.  
We showed (Theorem 6.1 in \cite{Ban9}) that if {\bf K} is a co-elementary
class that is {\it co-inductive\/}, i.e., closed under limits of
0-chains, and if $X \in \mbox{\bf K}$ is infinite, then there is a compactum
$Y \in \mbox{\bf K}^1$, of the same weight as $X$, such that $X$ is a
continuous image of $Y$. (So $\mbox{\bf K}^1$ is quite substantial under
these circumstances.)  
{\bf CH}, {\bf BS}, and {\bf CON} (the class of 
{\it continua\/}, i.e., connected compacta) are easily seen to be examples
of co-inductive co-elementary classes.  In \cite{Ban9} 
we showed $\mbox{\bf CH}^1 = \mbox{\bf BS}^1 = \{\mbox{Boolean spaces
without isolated points}\}$
(Proposition 6.2), and that every member of $\mbox{\bf CON}^1$
(i.e., every {\it co-existentially closed continuum\/})
is {\it indecomposable}, i.e., incapable of being written as
the union of two proper subcontinua (Proposition 6.3).  We posed the
question of whether $\mbox{\bf CON}^1$  is a co-elementary class,
and conjectured
that every co-existentially closed continuum is of (Lebesgue covering)
dimension one.  While the question of co-elementarity is still open, we
have been able to settle the conjecture in the affirmative.  We are
grateful to Wayne Lewis \cite{Lew}, who suggested the use of a theorem
of D. C. Wilson \cite{Wils}.\\
  
\subsection{Theorem.}\label{4.5}
Every co-existentially closed continuum is an indecomposable continuum
of dimension one.\\

\noindent
{\bf Proof.}  Because of Proposition 6.3 of \cite{Ban9}, we need only
concentrate on the issue of dimension.

Let $Q$ denote the Hilbert cube, the usual topological product of countably
many copies of the closed unit interval.  It is well known \cite{Wil} that  
every metrizable compactum can be replicated as a (closed) subspace of $Q$.
Next, let $M$ denote the Menger universal curve, a one-dimensional Peano
(i.e., locally connected metrizable) continuum.  Perhaps less well known
is the fact \cite{Nad} that every one-dimensional metrizable compactum
can be replicated as a (closed) subspace of $M$.  Wilson's
theorem \cite{Wils} says that there is a continuous surjection 
$f:M \to Q$ whose point-inverses are all homeomorphic to $M$.  So $f$ is,
in particular, monotone; hence inverse images of subcontinua of $Q$ are
subcontinua of $M$.  
Now let $X$ be any metrizable continuum, viewed
as a subspace of $Q$.  Then $f^{-1}[X]$ is a subcontinuum of $M$ that
maps via $f$ onto $X$.  Since $M$ is one-dimensional, so is $f^{-1}[X]$.

So we know that every metrizable continuum is a continuous image of a
metrizable continuum that is one-dimensional.  Let $X$ now be an arbitrary
continuum.  Then, by L\"{o}wenheim-Skolem, there is a co-elementary
map $f:X \to Y$, where $Y$ is a metrizable continuum.  Using the
result in the preceding paragraph, let $g:Z \to Y$ be a continuous surjection,
where $Z$ is a metrizable continuum of dimension one.  Because of
the co-elementarity of $f$, there is a homeomorphism 
$h:XI\backslash {\cal D} \to YI\backslash {\cal D}$ of ultracopowers
such that $f\circ p = q \circ h$, where $p$ and $q$ are the obvious
codiagonal maps.  Since covering dimension is an invariant of co-elementary
equivalence \cite{Ban2}, we know that $ZI\backslash {\cal D}$ is a
continuum of dimension one.  Thus $p\circ h^{-1}\circ (gI\backslash {\cal D})$
is a continuous surjection from a continuum of dimension one onto $X$.

Now suppose $X$ is 1-closed in {\bf CON}.  Then, by the paragraph above,
there is a continuous surjection $f: Y \to X$, where $Y$ is a continuum of
dimension one.  But $f$ is a co-existential map, and co-existential maps
preserve being infinite, and cannot raise dimension.  The dimension of
$X$ cannot be zero; hence it must be one. $\dashv$\\  

The following result records some general information concerning levels
of maps between classes, and is an easy corollary of the general
results above.\\

\subsection{Corollary.}\label{4.6}
Let {\bf K} be a class of compacta, $\alpha$ an ordinal.

$(i)$ Suppose $\alpha > 0$, and $\mbox{\bf K}^{\alpha}$ is closed under
ultracopowers.  Then
$\mbox{Lev}_{\geq 0}(\mbox{\bf K}^{\alpha},\mbox{\bf K}^{\alpha})= 
\mbox{Lev}_{\geq \omega}(\mbox{\bf K}^{\alpha},\mbox{\bf K}^{\alpha})$.

$(ii)$ Suppose {\bf K} is closed under ultracopowers.  Then 
$\mbox{Lev}_{\geq \alpha}(\mbox{\bf K}^{\alpha},\mbox{\bf K})= 
\mbox{Lev}_{\geq \alpha}(\mbox{\bf K}^{\alpha},\mbox{\bf K}^{\alpha})$.

$(iii)$ Suppose {\bf K} is closed under ultracopowers, and
$\mbox{Lev}_{\geq \alpha}(\mbox{\bf K},\mbox{\bf CH})= 
\mbox{Lev}_{\geq \alpha}(\mbox{\bf K},\mbox{\bf K})$.  Then
$\mbox{Lev}_{\geq \alpha +1}(\mbox{\bf K}',\mbox{\bf CH})= 
\mbox{Lev}_{\geq \alpha +1}(\mbox{\bf K}',\mbox{\bf K}')$ (where
$\mbox{\bf K}' := \mbox{\bf CH}\setminus \mbox{\bf K}$).\\

\noindent
{\bf Proof.} {\it Ad\/} $(i)$: By definition of $\mbox{\bf K}^{\alpha}$,
$\mbox{Lev}_{\geq 0}(\mbox{\bf K},\mbox{\bf K}^{\alpha})= 
\mbox{Lev}_{\geq \alpha}(\mbox{\bf K},\mbox{\bf K}^{\alpha})$.  The 
conclusion is immediate, by \ref{4.1}.

{\it Ad\/} $(ii)$: There is nothing to prove if $\alpha = 0$.  So assume
$\alpha >0$, and suppose $X$ is $\alpha$-closed in {\bf K}, $Y \in
\mbox{\bf K}$, and $f:X \to Y$ is of level $\geq \alpha$.  Let 
$p: YI\backslash {\cal D} \to Y$ and $g:YI\backslash {\cal D} \to X$
witness the fact; i.e., $p$ is a codiagonal map, $g$ is of level
$\geq \alpha -1$, and $f\circ g = p$.  Let $Z \in \mbox{\bf K}$, with
$h: Z \to Y$ a continuous surjection.  Let $q: ZI\backslash {\cal D} \to
Z$ be the appropriate codiagonal map.  Since {\bf K} is closed under
ultracopowers, and $X \in \mbox{\bf K}^{\alpha}$, we know that
both $g$ and $g\circ (hI\backslash {\cal D})$ are of level $\geq \alpha$.
Then $f\circ g\circ (hI\backslash {\cal D}) = p\circ (hI\backslash {\cal D})
= h\circ q$ is of level $\geq \alpha$, by \ref{2.5}.  By \ref{2.6}, 
$h$ is also of level $\geq \alpha$; hence $Y \in \mbox{\bf K}^{\alpha}$.   

{\it Ad\/} $(iii)$: Suppose $f:X \to Y$ is a map of level $\geq \alpha +1$,
and $Y \in \mbox{\bf K}$.  We need to show $X \in \mbox{\bf K}$.  But
this is immediate from the definition of level, plus our hypotheses. $\dashv$\\

\subsection{Remark.}\label{4.7} We have very few results concerning the
nature of $\mbox{\bf K}^{\alpha}$, given information about {\bf K}.
We can prove quite easily, though, that $\mbox{\bf CH}^2$, 
$\mbox{\bf BS}^2$, and 
$\mbox{\bf CON}^2$ are all empty.  Indeed, let $X$ be any compactum, with
$Y$ the disjoint union of $X$ with a singleton, and $Z$ the product of
$X$ with a Cantor discontinuum.  Then there exist continuous surjections
$f:Y \to X$ and $g:Z \to X$.  Assume $X$ is now 2-closed in {\bf CH}.  
Maps of level $\geq 1$ preserve the
property of having no isolated points (Proposition 2.8 in \cite{Ban9});
so we conclude that $X$ has
no isolated points because $Z$ has none.  On the other hand, since the
class of compacta without isolated points is co-elementary, and
$f$ is of level $\geq 2$,  
we conclude, by \ref{4.6}$(iii)$,
that $X$ has an isolated point because $Y$ does.  Thus 
$\mbox{\bf CH}^2$ is empty.  If $X$ above happens to be Boolean, so are
$Y$ and $Z$; hence the same argument shows that 
$\mbox{\bf BS}^2$ is empty.  Now assume $X \in \mbox{\bf CON}^2$.  Then
$X$ has dimension one, by \ref{4.5}.  Let $Y$ be the product of $X$
with the Hilbert cube.  Then there is a continuous surjection
$f:Y \to X$, and $Y$ is an infinite-dimensional continuum.  Since the
class of finite-dimensional continua is co-elementary, as well as closed
under images of maps of level $\geq 1$ (Proposition 2.6 in \cite{Ban9}), 
and $f$ is of
level $\geq 2$, we conclude, again by \ref{4.6}$(iii)$, that
$X$ is infinite-dimensional because $Y$ is.  Thus   
$\mbox{\bf CON}^2$ is empty.\\

\end{document}